\documentclass[11pt]{article}
\usepackage[T1]{fontenc}
\usepackage{times}
\usepackage{mathdots}
\usepackage{multirow}
\usepackage{textcomp}
\usepackage{amssymb,amscd,amsthm,latexsym}
\usepackage{amsmath}
\usepackage{graphicx}
\usepackage{float}
\usepackage{amssymb}
\usepackage{epstopdf}
\usepackage{verbatim}

\newtheorem{theorem}{Theorem}[section]

\newtheorem{lemma}[theorem]{Lemma}

\theoremstyle{definition}
\newtheorem{definition}[theorem]{Definition}

\newtheorem{example}[theorem]{Example}

\newtheorem{cor}[theorem]{Corollary}

\newtheorem{rem}[theorem]{Remark}

\def\n{n \times n}

\def\b{\beta}
\def\l{\lambda}

\def\s{\sigma}
\def\e{\epsilon}
\title{\vskip10mm An application of the geometry of complex polynomials in the theory of location of matrix eigenvalues
}
\author{Frank J. Hall\\ Department of Mathematics and Statistics\\ Georgia State University\\
	Atlanta, GA, USA\\ fhall@gsu.edu
	\, \\~~~\\
	Rachid Marsli
	\\ Department of Mathematics\\ Perimeter College \\ Georgia State University\\ Atlanta, GA, USA\\
	rmarsli1@gsu.edu
	\, \\~~~\\
	Michael Stewart\\ Department of Mathematics and Statistics\\ Georgia State University\\
	Atlanta, GA, USA\\ mastewart@gsu.edu
}
\date{\today}
\begin{document}
	%
	%
	%
	\maketitle
	\begin{abstract}
		\noindent
		Let $A$ be an $n \times n$ complex matrix with $k$ known eigenvalues. In this work, it is shown how to use those $k$ eigenvalues to obtain a localization region of the remaining $(n-k)$ eigenvalues of $A$. As a consequence, a clustering theorem is obtained. Other related ideas are discussed. Small examples and illustrations are given to illustrate the general principals.  A larger symmetric example shows that the method can be used to identify larger features of the spectrum, such as spectral gaps, using a relatively small number of eigenvalues computed by the Lanczos algorithm.
	\end{abstract}
	\noindent
	{\it AMS Subj. Class.:} 15A18, 30C10, 65F99
	\section{Introduction}
		Localizing the spectrum of real and complex matrices within relatively small regions of the complex plane is a fundamental problem in matrix analysis. Among the most well-known results in this area are the Gershgorin disks, which provide circular regions centered at the diagonal entries of a matrix \cite{G, Horn, V}, and the Brauer sets, which consider pairs of distinct rows and yield tighter inclusion regions in the form of Cassini ovals \cite{Brauer, Horn, V}.
		
		More recently, Frank J. Hall and Rachid Marsli have developed a sequence of results extending the classical Gershgorin theory. In \cite{M1}, they proved that if a complex matrix has an eigenvalue of geometric multiplicity $k$, then that eigenvalue belongs to at least $k$ Gershgorin disks. In \cite{M3}, they introduced the concept of Gershgorin disks of the second type. Subsequently, they proved in \cite{Hall1} that every nontrivial eigenvalue of a real matrix $A$ with constant row sums or constant column sums is contained in a Gershgorin disk of the second type of  the transpose of $A$.
		
		It is difficult to provide a comprehensive account of the many interesting contributions to the theory of location of eigenvalues. Nonetheless, we mention here several other representative articles \cite{Barany, Dall, De, Hall2, K, M2, M4, Me, Pena, W}. Two recent surveys providing overviews of the development of this theory are \cite{Rehman, Z}. Two valuable books are \cite[Chapter 6]{Horn} and \cite{V}.
		
		In \cite{M2} and \cite{Me}, the authors discuss how a given real eigenvalue, together with its eigenvectors, can be used to obtain a localization region for the spectrum. The present work follows the same general spirit, but approaches the problem from a fundamentally different perspective. We exploit the geometry of complex polynomials to address the following question:\\	
		Let $A=[a_{ij}]$ be an $n\times n$ complex matrix. Suppose that $k$ eigenvalues of $A$ are known. How can this information be used to localize the remaining $n-k$ eigenvalues of $A$?
	\section{An inclusion region of the spectrum of a matrix with some known eigenvalues}
	If some of the eigenvalues of a complex matrix $A$ are known, then the rest of its spectrum can be located using the following theorem.
	\begin{theorem}\label{t1}
		Let $A$ be an $\n$ nonzero complex matrix for some integer $n\geq3$ with spectrum $\l_1, \dots, \l_n$ counting multiplicities.
		For some positive integer $k<n$, suppose that $\l_1, \dots, \l_k$ are distinct. Let $B$ be the matrix given by 
		\begin{equation}\label{eq4}
			B = \prod_{j=1}^{k} (A-\l_j I),
		\end{equation}
		where $I$ is the $\n$ identity matrix and let $\rho(B)$ be the spectral radius of $B$.
		Then all the eigenvalues of $A$ lie within the region of the complex plane defined by the inequality
		\begin{equation}\label{eq5}
			\prod_{j=1}^{k}|z-\l_j| \leq \rho(B),
		\end{equation}
		where $z$ is the complex variable.
	\end{theorem}
	\begin{proof}
		Let $A$ be an $\n$ nonzero complex matrix for some integer $n\geq3$ with spectrum $\l_1, \dots, \l_n$ counting multiplicities.
		For some positive integer $k<n$, suppose that $\l_1, \dots, \l_k$ are distinct. Let $B$ the matrix given by 
		$$
			B = \prod_{j=1}^{k} (A-\l_j I).
		$$
		Then it can be seen that the eigenvalues of $B$ are given by $\b_1 = \dots = \b_k=0$ and $\b_i = \prod_{j=1}^{k}(\l_i - \l_j)$ for $i=k+1, \dots, n$.
		This important fact holds whether $A$ is diagonalizable or defective. Then for every eigenvalue $\beta_i$ of $B$, we have
	 	\begin{equation}\label{eq3}
	 		\prod_{j=1}^{k}|\l_i - \l_j| \leq \rho(B), ~~ 1\leq i \leq n.
	 	\end{equation}	 	 
	 	Therefore, $\l_i$ is in the region defined by the inequality in (\ref{eq5}).	 	
	 \end{proof}
	\noindent
	Now consider the polynomial $p(z) = \prod_{j=1}^{k}(z-\l_j)$. From the inequalities in (\ref{eq5}), we see that the eigenvalues of $A$
	lie within the region of the complex plane defined by the inequality $|p(z)|\leq \rho(B)$. This type of region has been well studied in both classical and modern geometry. Its boundary is known as a polynomial lemniscate or a generalized Cassini oval.\\
	
	\noindent
	\begin{rem}\label{r1} ~
			\begin{enumerate}
			\item If $k=1$ then the inequality in (\ref{eq5}) reduces to the disk with center at $\l_1$ and radius $\rho(B)$.
			\item If $k=2$ then the inequality in (\ref{eq5}) becomes $|z-\l_1| |z-\l_2| \leq \rho(B)$. This defines the region well-known as a Cassini oval. However, what we are doing here is different from the Cassini ovals related to Brauer sets, which is another context in which location of eigenvalues is studied in \cite{Brauer} and \cite[Chapter 6]{Horn}. 
			\item As $k$ increases, the region becomes more complicated and difficult to predict, although it has in general some very nice geometric properties due to the holomorphic aspect of complex polynomials.
			\item As $\rho(B)$ varies, the shape of the region changes significantly. For very small $\rho(B)$, the region consists of several small, disjoint ovals each of which surrounds one $\l_j$ for $j=1, \dots, k$. As $\rho(B)$ increases, the ovals increase in size and start to merge one after another into one area that, when $\rho(B)$ is large enough, looks like a disk.
			\item The region defined by (\ref{eq5}) cannot have more than $k$ pairwise disjoint ovals. This is a well-known fact in complex polynomial theory and it is a consequence of the fundamental theorem of algebra and the maximum modulus principle.
			\item According to the maximum modulus principle, a number $z$ for which $|p(z)|=\rho(B)$ must lie on the boundary of this region.
			\item It follows from the above remark that there is always an eigenvalue of $A$ on the boundary of the region given by Theorem \ref{t1}.
			In fact, there is always at least one eigenvalue $\l_i$ of $A$ for which  $\prod_{j=1}^{k}|\l_i-\l_j| = \rho(B)$. Therefore, by the maximum modulus principle, such eigenvalue $\l_i$ must lie on the Boundary of the region.
		\end{enumerate}
	\end{rem}
	\begin{example}\label{ex1}
		Consider the $8 \times 8$ matrix
		$$	
			A\,=\, \left[
		\begin{array}{r r r r r r r r }
			17 & -1 &   -15 &  -2  & 0 &   -15 &  -39  & 63\\
			0 & -4 &   1 &  -12 & 0 &   1 &  1  & -1\\
			-7 & 3 &  22.5  &  6  & 2 &   17 &  35.5  &  -42.5\\
			0 & 6 &   -1 &  14  & 0 &   -1 &  -1  &   1\\
			0 & 1 &   1.5 &  2 &  11.5 &   6.5 &  1.5&   -1.5\\
			0 & -1 &   -0.5 &  -2 &  -2 &  5  &  -0.5 &   0.5\\
			-2 & -1 &   -13 &  -2 &  0 &   -13 &  -16 & 7\\
			-9 & 0 &   1 &  0 &  0 &   1 &  11 &   -27
		\end{array}
		\right].
		$$
		\noindent
		The eigenvalues of this matrix are: $1 , 2, -7, -6, 8\pm i$ and $8.5 \pm i$.
		If none of these eigenvalues is known we can locate them using a disk centered at the origin with a radius equal to one of the norms of $A$, We can also use the famous Gershgorin theorem, according to which all the eigenvalues are included in the union of the disks with centers and radii, respectively, equal to the diagonal element in each row and the sum of the absolute values of the off-diagonal entries on the same row. By a quick look at the matrix, we notice the radii of the disks obtained by using either a norm or the Gershgorin theorem are relatively large compared to the size of the eigenvalues. Other techniques, such as the one explained in \cite{Hall2} and based on Gelfand's formula can be used if a good upper approximation of the dominant eigenvalue(s) is known.
		If any eigenvalue of $A$ is known, say $\l = 2$, then according to our theorem, all the other eigenvalues are located within the disk defined by the inequality
		$$ |z-2| \leq \rho(B_1),$$ where $\rho(B_1)$ is the spectral radius of $B_1 = A - 2I$. This is shown in the figure below.
		\begin{figure}[H]
			\centering
			\includegraphics[scale=0.30]{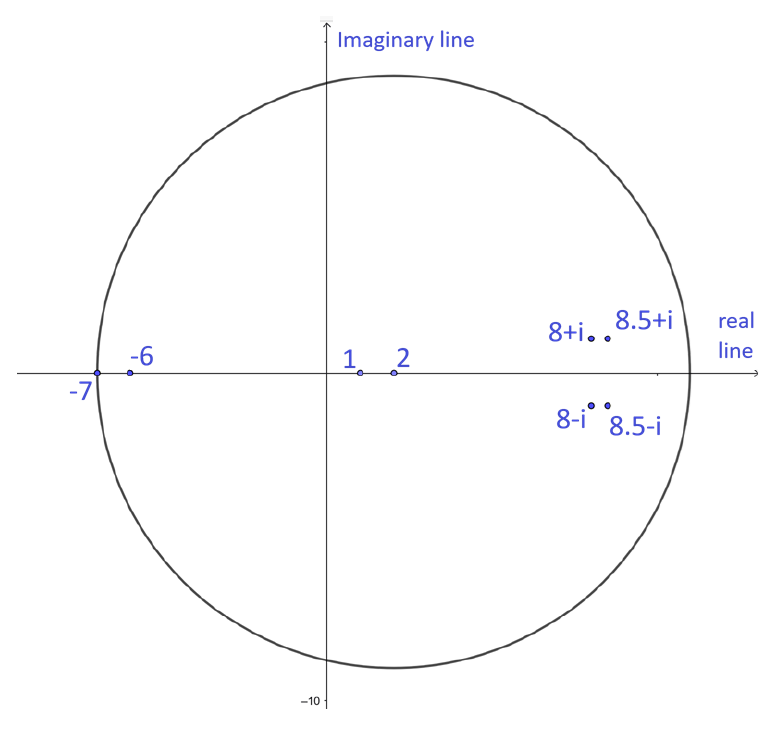}			
			{\scriptsize
				\begin{minipage}{0.7\textwidth}
					\caption{}
					\vspace{2mm}
					\centering
					All the eigenvalues of $A$ lie within the disk with center at $\l =2$ and radius $\rho(B_1) = 9$.
				\end{minipage}
			}					
		\end{figure}
		\noindent
		Now if two distinct eigenvalues of $A$ are known, let's say $\l_1 = -7$ and $\l_2 = 2$, then all the remaining eigenvalues of $A$ are located within the region defined by the inequality 
		$$
		|z+7| |z-2| \leq \rho(B_2),
		$$
		where $\rho(B_2)$ is the spectral radius of the matrix $B_2 = (A+7I)(A-2I)$. Note that all the eigenvalues of $B_2$ have the form $(\l_i + 7)(\l_i-2)$ where $\l_i$ ranges over all the eigenvalues of $A$. This allows easy calculation of the spectral radius of $B_2$ which is $\rho(B_2) \approx 102.15$. Below is the obtained inclusion region of the eigenvalues of $A$.
				\begin{figure}[H]
					\centering
					\includegraphics[scale=0.30]{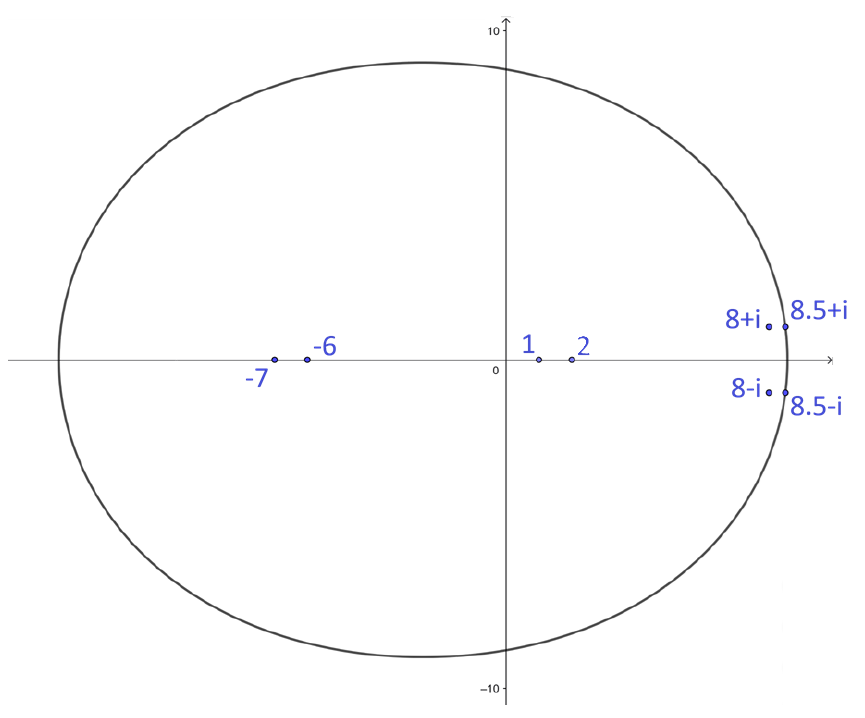}			
					{\scriptsize
						\begin{minipage}{0.7\textwidth}
							\caption{}
							\vspace{2mm}
							\centering
							All the eigenvalues of $A$ lie within the area defined by $|z+7| |z-2| \leq \rho(B_2) \approx 102.15$.
						\end{minipage}
					}					
				\end{figure}
				\noindent
				This region seems close to the first one in size but a little larger. One may expect a smaller region since we are using more information; two eigenvalues instead of one. Actually, the size of the obtained region depends on many factors including the entries of the relative positions of the known eigenvalues, how far they are from each other and how large or small is the spectral radius of the constructed matrix, $B_2$ in this case. Choosing another pair of eigenvalues, for example $\l_1 = -7$ and $\l_2 = 8+i$, results in a neatly smaller region defined by
				$$
				|z+7|~|z-8-i| \leq \rho(B_3),
				$$
				where $\rho(B_3)$ is the spectral radius of $B_3 = (A+7I)(A-(8+i)I)$. See the obtained region below.
						\begin{figure}[H]
						\centering
						\includegraphics[scale=0.40]{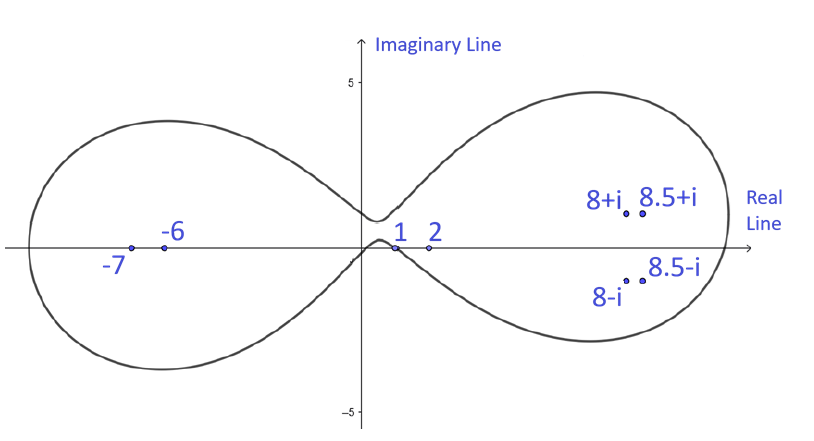}			
						{\scriptsize
							\begin{minipage}{0.7\textwidth}
								\caption{}
								\vspace{2mm}
								\centering
								All the eigenvalues of $A$ lie within the area defined by $|z+7|~|z-8-i| \leq \rho(B_3) \approx 56.57$.
							\end{minipage}
						}					
					\end{figure}		
				\noindent
				Suppose that we know three eigenvalues of $A$, $\l_1 = -7, \l_2 = 2$ and $\l_3 = 8+i$. Then all the other eigenvalues lie within the region defined by $$|z+7|~|z-2|~|z-8-i| \leq \rho(B_4),$$ where $B_4 = (A+7I)(A-2I)(A-(8+i)I)$ and $\rho(B_4) = 210.58$.
					\begin{figure}[H]
					\centering
					\includegraphics[scale=0.40]{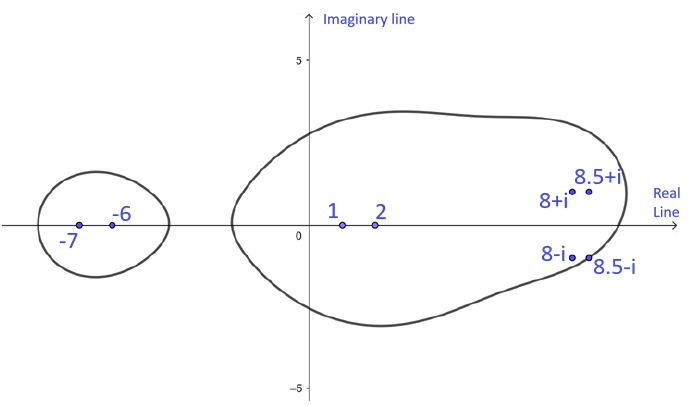}			
					{\scriptsize
						\begin{minipage}{0.7\textwidth}
							\caption{}
							\vspace{2mm}
							\centering
							All the eigenvalues of $A$ lie within the area defined by $|z+7|~|z-2|~|z-8-i| \leq \rho(B_4) \approx 210.58$. The region is made of two disjoint areas.
						\end{minipage}
					}					
				\end{figure}						
				\noindent
				Note that there is at least one eigenvalue of $A$ on the boundary of each of the obtained regions. This is in accordance with Remark \ref{r1}\,(7).
	\end{example}
	\section{Application to matrices with real spectrum}
		Many real and complex matrices possess purely real spectra, including well-known classes such as Hermitian, symmetric and positive definite matrices.
		\begin{cor}\label{c2}
				Let $A$ be an $\n$ nonzero complex matrix for some integer $n\geq3$ with spectrum $\l_1, \dots, \l_n$ counting multiplicities.
			For some positive integer $k<n$, suppose that $\l_1, \dots, \l_k$ are distinct and the remaining eigenvalues $\l_{k+1}, \dots, \l_n$ are real.
			Let $B$ the matrix given by 
			\begin{equation}\label{eq4}
				B = \prod_{j=1}^{k} (A-\l_j I),
			\end{equation}
			where $I$ is the $\n$ identity matrix and let $\rho(B)$ be the spectral radius of $B$.
			Then $\l_{k+1}, \dots, \l_n$ lie on the interval or union of intervals resulting from the intersection of the real line and the region of the complex plane defined by the inequality
			\begin{equation}\label{eq}
				\prod_{j=1}^{k}|z-\l_j| \leq \rho(B),
			\end{equation}
			where $z$ is a complex variable. 
		\end{cor}
	\begin{example}\label{ex3}
		Consider the $8 \times 8$ symmetric matrix
		$$	
		A\,=\, \left[
		\begin{array}{r r r r r r r r }
			1 & 5 &   -5 &  7  & 1 &   -3 &  3  & -1\\
			5 & -7 &   -5 &  7 & 5 &   1 &  3  & -1\\
			-5 & -5 &  1  &  9  & 11 &   3 &  -7  &  1\\
			7 & 7 &   9 &  1  & -9 &   -1 &  1&   -7\\
			1 & 5 &   11 &  -9 &  -7 &   5 &  -5&   7\\
			-3 & 1 &   3 &  -1 &  5 &  1  &  3 &   -1\\
			3 & 3 &   -7 &  1 &  -5 &   3 &  9 & 1\\
			-1 & -1 &   1 &  -7 &  7 &   -1 &  1 &   9
		\end{array}
		\right].
		$$		
		\noindent
		The eigenvalues of $A$, rounded to the second decimal digit, are given by $-27.15$, $-9.01$, $-8$, $0$, $6.06$, $8$, $16$, and $22.11$.
		Suppose that the smallest and largest eigenvalues $\l_1 \approx -27.15$ and $\l_2 \approx 22.11$ are known to us and we want to use this information to locate the eigenvalues of $A$. We construct the matrix $B = (A+27.15I)(A-22.11I)$, which has a spectral radius $\rho(B) \approx 600.29$.
		Then according to Theorem \ref{t1}, the remaining eigenvalues of $A$ lie within the region defined by the inequality
		$$
		|z+27.15| |z-22.11| \leq 600.29.
		$$
		See Figure $5$ below.
		\noindent
		\begin{figure}[H]
			\centering
			\includegraphics[scale=0.40]{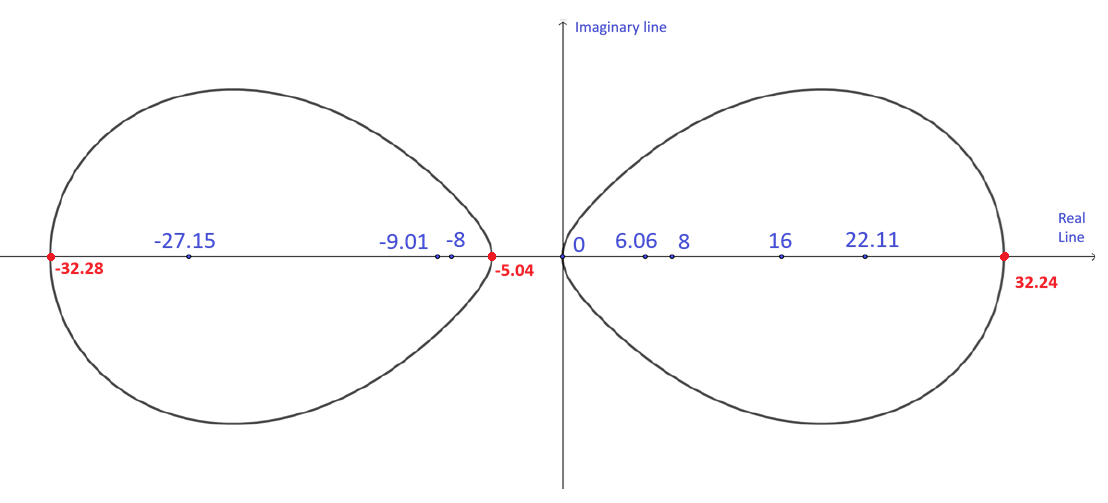}			
			{\scriptsize
				\begin{minipage}{0.7\textwidth}
					\caption{}
					\vspace{2mm}
					\centering
					All the eigenvalues of $A$ lie within the area defined by $|z+27.15| |z-22.11| \leq 600.29$. The region is made of two disjoint areas.
				\end{minipage}
			}					
		\end{figure}						
		\noindent
		Since $A$ is symmetric, all its spectrum is real and Corollary \ref{c2} ensures that all the remaining eigenvalues lie within the union of intervals
		$$[-37.28 \, , \, -5.04] \cup [0\, , \,32.24],$$
		as in the figure below.	Note that one of the boundary points of this union of intervals is the eigenvalue $\l = 0$ of $A$, while the other boundaries $-32.28, -5.04$ and $32.24$ are geometrically determined.
		\begin{figure}[H]
			\centering
			\includegraphics[scale=0.40]{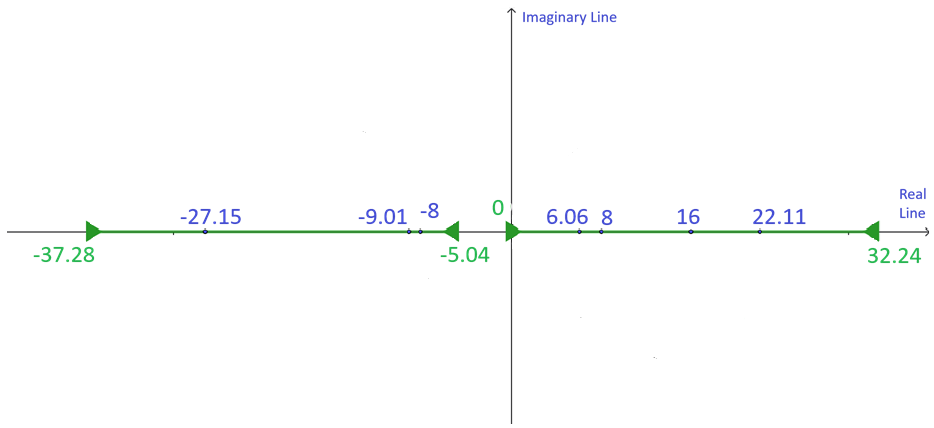}			
			{\scriptsize
				\begin{minipage}{0.8\textwidth}
					\caption{}
					\vspace{2mm}
					\centering
					All the eigenvalues of $A$ lie within the union of intervals\\ $[-37.28 \, , \, -5.04] \cup [0\, , \,32.24]$.					
				\end{minipage}
			}					
		\end{figure}						
		\noindent
	\end{example}
	\begin{rem}
		It is important to note that at least one of the boundaries of the interval or union of intervals given by Corollary \ref{c2} is an eigenvalue of the original matrix $A$ (Example \ref{ex3}, $\l =0$). Theorem \ref{t1} or Corollary \ref{c2}, by themselves do not provide ways to determine the values of the other boundaries. For that, one can leverage the geometric properties of the obtained region as we did in Example \ref{ex3}, or if possible, the algebraic or numerical properties of the corresponding inequality.
	\end{rem}
	\begin{rem}
		If $\l_1, \dots, \l_k$ are real, then the region of the complex plane defined by
		$|z-\l_1| |z-\l_2| \dots |z-\l_k| \leq \rho(B)$ exhibits a symmetry with respect to the real axis. In fact, if a point $M_1(x,y)$ belongs to this region, then clearly the point $M_2(x,-y)$ also lies in the same region.	This explains why in the previous  example, the obtained region has symmetry with respect to the real axis.
	\end{rem}
	\section{An eigenvalue clustering theorem}
		The expression “clusters of eigenvalues” generally refers to the distribution of eigenvalues of a given matrix in the complex plane and, in particular, to their relative proximity to one another. Sometimes, eigenvalues appear to gather into several groups, such that the eigenvalues within each cluster are relatively close to one another, while those belonging to different clusters are relatively far apart. This concept has been studied and applied in numerical linear algebra \cite{Chan, CS, Chu, S} as well as in statistical data analysis \cite{Han, L, Tao}. A formal definition may vary slightly from one reference to another. In our current work, we adopt the following definition.
	\begin{definition}
		Let $A$ be an $\n$ nonzero complex matrix with spectrum $\{\l_1, \dots,$ $\l_n\}$ counting multiplicities.
		A clustering set of the eigenvalues of $A$ or a spectral cluster set of $A$ is a union of disjoint, isolating regions in the complex plane that partition the spectrum of $A$.\\
		In particular, for some integer $k <n$, the interior of the set $\{D_i(\s_i, r_i) ~~|~~ 1\leq i \leq k\}$ is a clustering set of the eigenvalues of $A$ if the disks are pairwise disjoint except for tangencies, every disk contains at least one eigenvalue, and every eigenvalue lies in the interior of a disk.
	\end{definition}
	\begin{definition}\label{d2} \cite[Section 6.1]{PR}
		Let $\{\l_1, \dots, \l_k\}$ be a set of $k$ points in the complex plane.
		A tree connecting the points $\l_1, \dots, \l_k$ is said to be a  Euclidean minimum spanning tree (EMST) if it uses the shortest possible total edge weight, where the weight of each edge is the Euclidean distance between its endpoints. 
	\end{definition}
	\noindent
	An EMST has many well-known geometric and combinatorial (graph theoretic) properties.
	If an EMST is embedded in the complex plane (which is Euclidean), then every vertex is connected to at least one of the closest vertices. Moreover, the circle whose diameter is an edge of the tree contains no vertex other than the endpoints of that diameter. Another well-known geometric property follows consequently from the previous two and is stated as a lemma for its use in proving the subsequent theorem.
	\begin{lemma}\label{l1}
		Let $\mathbf{T}$ be a  Euclidean minimum spanning tree connecting a set of $k$ points $\{\l_1, \dots, \l_k\}$ in the complex plane,
		connected by $(k-1)$ edges $d_1, \dots, d_{k-1}$, where $k\geq3$.
		For $j=1, \dots, k$, define $r_j$ to be half of the minimum length among the edges incident to $\l_j$, that is
		\begin{equation}\label{eq2}
			r_j = \frac{1}{2} \min_{1\leq i < k} \left\{ |d_i| ~~\big|~~ d_i \text{~~is incident to~~} \l_j \right\},
		\end{equation}			
		where $|d_i|$ is the length (weight) of the edge $d_i$.
		Then the interiors of the disks $D_j(\l_j, r_j)$, for $j = 1, \dots, k$ are pairwise disjoint.
	\end{lemma}
	\begin{theorem}\label{t2}
		Let $A$ be an $n \times n$ nonzero complex matrix with $n \geq 3$.
		Let $\l_1, \dots, \l_n$ be the eigenvalues of $A$ counting multiplicities.
		For $ 2 \leq k <n$ (practically $k<<n$), suppose that $\l_1, \dots, \l_k$ are distinct 
		and connected by a Euclidean minimum spanning tree whose edges are $d_1, \dots, d_{k-1}$.		
		Consider the product
		\begin{equation}\label{eq10}
			r= \prod_{j=1}^{k} r_j
		\end{equation}	
		\noindent
		and the matrix 
		\begin{equation}\label{eq7}
			B = \prod_{j=1}^{k} (A - \l_j I),
		\end{equation}
		where $I$ is the $n \times n$ identity matrix and let
		$\rho(B)$ be the spectral radius of $B$.\\
		If 
		\begin{equation}\label{eq1}
			\rho (B) < r,
		\end{equation}
		then the set of disks $\{D_j(\l_j , r_j) \mid j=1, \dots, k \}$ form a clustering set of all the eigenvalues of $A$.
	\end{theorem}
	\begin{proof}				
		The eigenvalues of $B$ are given by $\b_1 = \dots = \b_k = 0$ and 
		$\b_i = \prod_{j=1}^{k} (\l_i - \l_j)$ 
		for $i=k+1, \dots, n$.
		Let $r = \prod_{j=1}^{k} r_j$ and
		suppose that there is some index $l \in \{k+1, \dots, n\}$ such that $\l_l$ does not lie in the interior of any disk $D_j$. Then $|\l_l - \l_j| \geq r_j$ for every $j \in \{1, \dots, k\}$ and consequently $\prod_{j=1}^{k} |(\l_l - \l_j)| \geq r$ or equivalently $|\b_l| \geq r$. But this implies that $\rho \geq r$ as $\rho(B)$ dominates all the eigenvalues of $B$.
		Therefore, if $\rho(B) < r$ then every eigenvalue must lie within the interior of a disk $D_j$.\\
		Moreover, by Lemma \ref{l1}, the interiors of those disks are pairwise disjoint, so they form a clustering set of the eigenvalues of $A$.
		To complete the proof we need to make evidence of the existence of matrices for which (\ref{eq1}) is satisfied. This can be done with a simple example. Consider a $3 \times 3$ matrix $A$ with eigenvalues $\l_1 =1, \l_2 =4$ and $\l_3=2$. Using the first two eigenvalues to locate the third one, we have $r_1 = r_2 = 1.5$, $B = (A-\l_1I)(A-\l_2I)$ and $\rho(B) =2 <r_1r_2 =2.25$.		
	\end{proof}
	\begin{rem}
		The converse of the assertion of Theorem \ref{t2} is not true in general. That is when all the eigenvalues lie within the interiors of the disks, the inequality $\rho(B) < r$  may or may not hold. 
	\end{rem}
	\begin{example}\label{ex2}
		Consider the $10 \times 10$ matrix
		$$	
		A\,=\, \left[
		\begin{array}{r r r r r r r r r r }
			-3  & 175 &   -263 &  162  & -385  &   250 &  -445  & 280 & -439  & 285 \\
			-18 & 190 &   -263 &  162  & -385  &   250 &  -445  & 280 & -439  & 285 \\
			-4  & 48  &   -74  &  50   & -94   &   50  &  -94   & 50  & -94   & 50  \\
			-4  & 58  &   -91  &  54   & -115  &   68  &  -122  & 68  & -122  & 68  \\
			-4  & 55  &   -82  &  52   & -124  &   87  &  -138  & 87  & -166  & 87  \\
			-4  & 51  &   -76  &  58   & -134  &   102 &  -157  & 119 & -218  & 119 \\
			-4  & 52  &   -81  &  66   & -147  &   125 &  -199  & 170 & -180  & 170 \\
			-4  & 39  &   -48  &  26   &  -74  &   52  &  -100  &  61 & -81   & 47 \\
			-4  & 16  &   -9  &  -6    & -3    &   -12 &   3    & -27 &  28   & -22  \\
			-2  & 4   &   -4   &  4    & -8    &    8  &  -12   & 12  & -16   & 31
		\end{array}
		\right]
		$$
		\noindent
		The eigenvalues of this matrix are: $15 ,~ 17,~ \pm 7i,~ \pm 10i,~ 16 \pm i$~ and ~$1 \pm 9i$.
		Suppose that we know only the three eigenvalues $\l_1 = 10i, \l_2 = -10i, \l_3 =17$ and we would like to locate the other eigenvalues of $A$. We first find a Euclidean minimum spanning tree connecting those three eigenvalues and the disks indicated by Theorem \ref{t2}, as shown in the Figure $7$.
		\begin{figure}[H]
			\centering
			\includegraphics[scale=0.50]{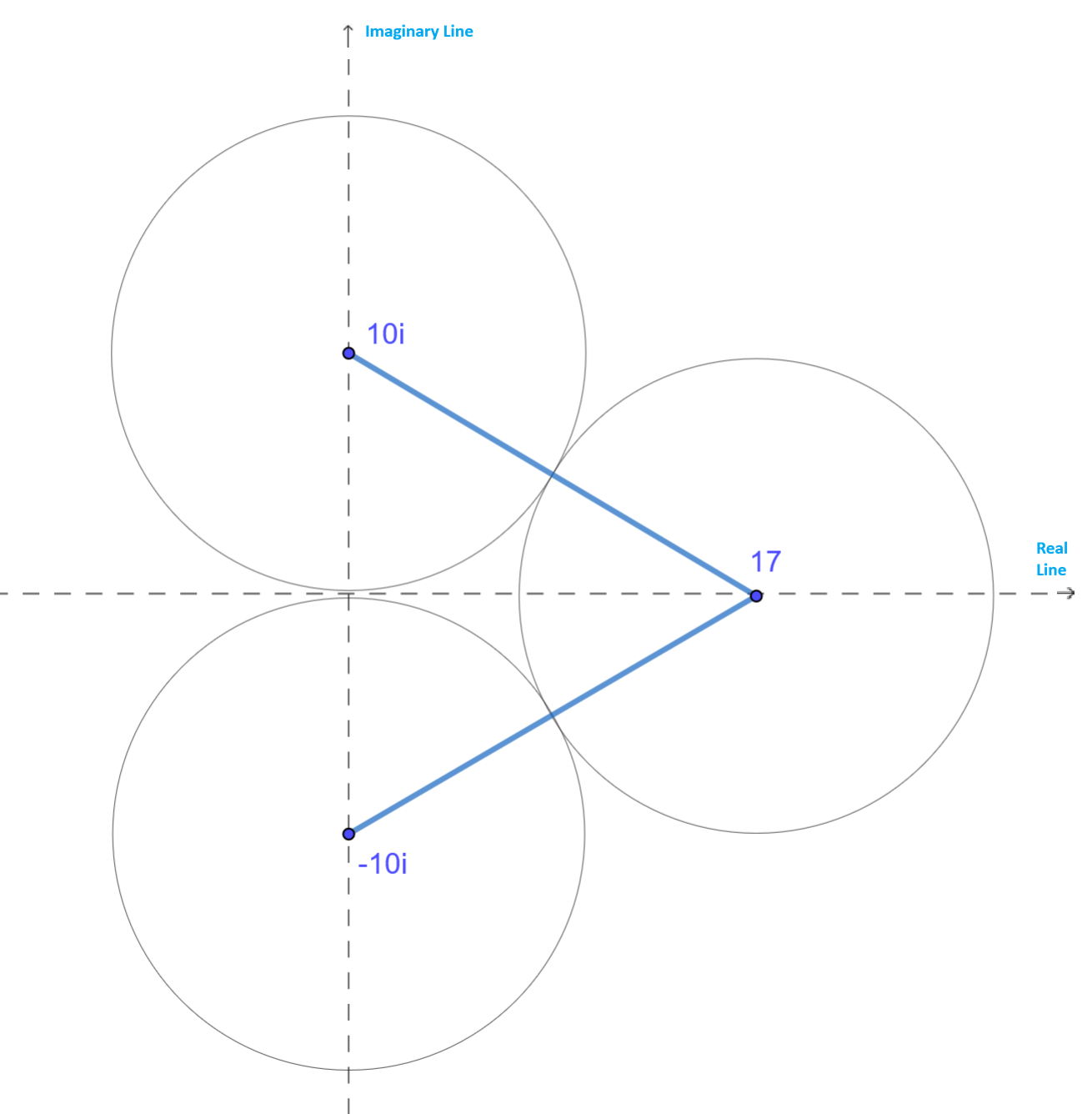}			
			{\scriptsize				
				\begin{minipage}{0.7\textwidth}
					\caption{}
					\vspace{2mm}
					\centering
					The line segments in blue form the  Euclidean minimum spanning tree connecting the three eigenvalues $\pm 10i$ and $17$.
					Each disk is centered at one of those eigenvalues and has a radius equal to half of the minimum edge length from its center to the eigenvalues connected to it by the spanning tree. In this special case, all the radii are approximately equal to $9.86$. 
				\end{minipage}
			}					
		\end{figure}	
		\noindent			
		A calculation of the product of the radii of the disks gives
		 $r_1 .r_2 . r_3 = \left(\frac{\sqrt{17^2 + 10^2}}{2}\right)^3\approx 959$.
		The matrix $B = (A-(10i)I)(A+(10i)I)(A-17I)$ has a spectral radius $\rho(B) \approx 937.62$.
		Since $\rho(B) < \prod_{i=1}^{3} r_i$, we conclude by Theorem \ref{t2} that the remaining eigenvalues of $A$ cluster in the union of the interiors of the three disks centered at $\pm 10i$ and $17$ with radii $r_1=r_2=r_3 \approx 9.86$, as shown in Figure $8$.
			\begin{figure}[H]
			\centering
			\includegraphics[scale=0.50]{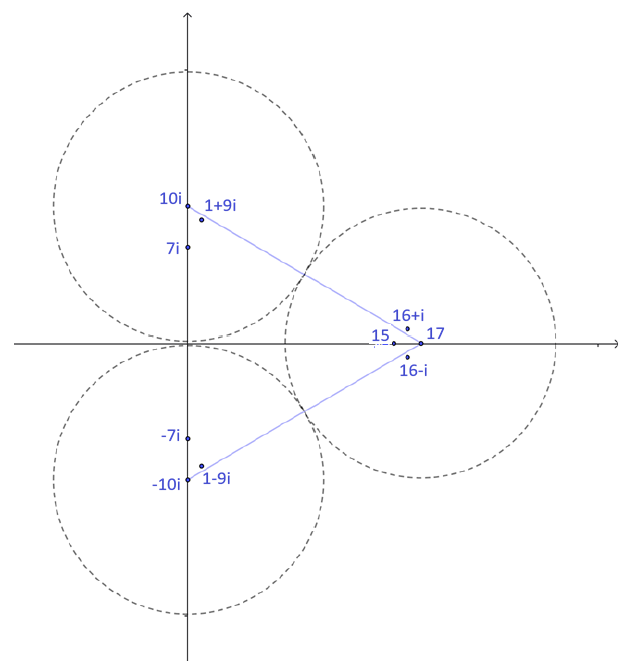}			
			{\scriptsize
				\begin{minipage}{0.7\textwidth}
					\caption{}
					\vspace{2mm}
					\centering
					All the eigenvalues of the matrix $A$ cluster in the interior of the disks. This is in accordance with Theorem \ref{t2}.
				\end{minipage}
			}					
		\end{figure}				
	\end{example}
	\begin{rem}\label{r2}
		Theorem \ref{t2} doesn't give a better location in size than the one given by Theorem \ref{t1}. It just provides a sufficient condition, $\rho(B) < r$, under which the spectrum of the matrix $A$ clusters within the interiors of disks $D(\l_j,r_j), j =1, \dots, k$, which have interiors disjoint from each other.
		When the condition $\rho(B)<r$ is satisfied, the Polynomial lemniscate defined by $\prod_{j=1}^{k} |z-\l_j| = \rho(B)$ subtends a region made up of $k$ disjoint parts known as generalized Cassini ovals, each of which 
		is centered at some $\l_j$ and contained within the interior of the disk $D(\l_j,r_j)$.		
	\end{rem}	
		\begin{cor}\label{c1}
		Let $A$ be an $n \times n$ nonzero complex matrix with $n \geq 3$.
		Let $\l_1, \dots, \l_n$ be the eigenvalues of $A$ counting multiplicities.
		For $ 2 \leq k <n$ (practically $k<<n$), suppose that $\l_1, \dots, \l_k$ are distinct. Let		
		\begin{equation}\label{eq6}
			\hat{r} = \frac{1}{2} \min \left\{ |\l_j - \l_i| ~~\big|~~ 1 \leq i < j \leq k \right\},
		\end{equation}			
		and consider the matrix $B$ given by (\ref{eq7}) and its spectral radius $\rho(B)$.
		If
		\begin{equation*}
			\rho(B) < \hat{r}^k,
		\end{equation*}
		then the set of disks $\{D_j(\l_j , \hat{r}) ~\mid~ j=1, \dots, k\}$ form a clustering set of all the eigenvalues of $A$.
	\end{cor}
	\begin{proof}
		From (\ref{eq2}) and (\ref{eq6}) we can see that $\hat{r} \leq r_j, j=1, \dots, k$. Moreover, Lemma \ref{l1} implies that the interior of the disks $\{D_j(\l_j , \hat{r})$ are disjoint. If $\rho(B) < \hat{r}^k$, then  by a reasoning similar to that in the proof of Theorem \ref{t2}, the set of disks $\{D_j(\l_j , \hat{r}) ~\mid~ j=1, \dots, k\}$ form a clustering set of all the eigenvalues of $A$.
	\end{proof}	
	\begin{rem}
		Corollary \ref{c1} is a consequence of Theorem \ref{t2}, but it doesn't supersede it. Before calculating the entire Euclidean minimum spanning tree, one can first compare  $\rho(B)$ to $\hat{r}^k$. If the condition $\rho(B) < \hat{r}^k$ is satisfied, then according to Corollary \ref{c1}, the spectrum of $A$ clusters within the interior of the disks $D(\l_j,\hat{r}), j=1, \dots, k$. Otherwise, one can  compare $\rho(B)$ to $r$, given by (\ref{eq10}).
		If the condition $\rho(B) < r$ is satisfied, then according to Theorem \ref{t2}, the spectrum of $A$ clusters within the disks $D(\l_j,r_j), j=1, \dots, k$.
	\end{rem}	
	\section{Further discussion of Theorem \ref{t2}}
	Going back to the matrix $A$ and Figure $8$ in Example \ref{ex2}, we can see that each of the three disks is significantly larger than where the included eigenvalues are located. 
	As in Example \ref{ex2}, suppose that we know only the three eigenvalues,  $\l_1 = 10i, \l_2 = -10i, \l_3 =17$ and we would like to locate the other eigenvalues of $A$, but this time by using Theorem \ref{t1} and compare the size of the obtained region to that of the disks obtained in Figure 8 by applying Theorem \ref{t2}. The region given by Theorem \ref{t1} consists of the three generalized Cassini ovals defined by
	$$
	|z-10i|~|z+10i|~|z-17| \leq \rho(B) \approx 937.62,
	$$
	as shown in the figure below. This area is neatly smaller than the union of the disks. One wonders whether this is just a coincidence or a mathematical fact is in play here and is behind this significant difference between the size of the disks and the generalized Cassini ovals.
	To answer this question, we recall two important ideas. The first is the well-known arithmetic mean geometric mean inequality (AM-GM inequality) which has applications in different areas of mathematics and the second is a result obtained by George P\'{o}lya in 1928 and provides an upper-bound on the area comprised within a polynomial lemniscate \cite{Poyla}.
	\begin{figure}[H]
		\centering
		\includegraphics[scale=0.85]{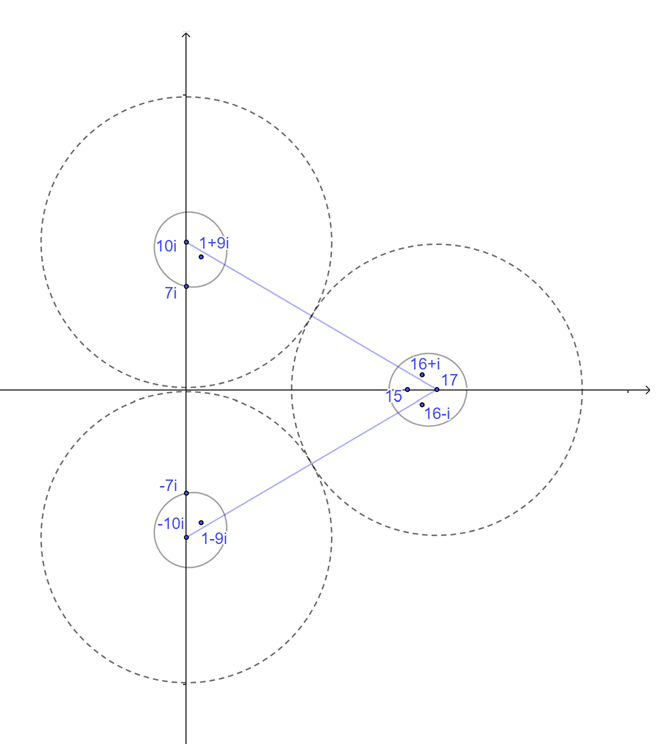}			
		{\scriptsize
			\begin{minipage}{0.7\textwidth}
				\caption{}
				\vspace{2mm}
				\centering
				Each generalized Cassini oval is contained within a disk, disjoint from the other ovals and contains a part of the spectrum.
				This is in accordance with Theorem \ref{t1}, Theorem \ref{t2} and Remark \ref{r2}. Moreover, the size of each oval is neatly smaller than 
				the disk in which it is contained.
			\end{minipage}
		}					
	\end{figure}				
	\begin{lemma}{(AM-GM inequality).}\label{l2}
		Let $a_1, \dots, a_k $ be positive numbers for some integer $k\geq 1$.
		Then
		$$ \frac{a_1+\dots+a_k}{k} ~~\geq~~ \sqrt[k]{\prod_{j=1}^k a_j}$$
	\end{lemma}
	\noindent
	\vskip 10pt
	Let $p(z)$ be a complex polynomial and $r$ be a positive number. The set of points 
	\begin{equation}\label{eq9}
		\varOmega = \left\{ z  \in \mathbb{C}~~ \big | ~~ |p(z)| \leq r \right\}
	\end{equation}
	is exactly the set of points contained within the polynomial lemniscate defined by $|p(z)|=r$ and has a well-known upper bound on its size, published by George P\'{o}lya in 1928, \cite{Poyla}.
	\begin{lemma}{(P\'{o}lya, 1928).}\label{l3}
		Let $p(z)$ be a complex polynomial having degree $k \geq 1$ and leading coefficient $a$. For any positive real number $r$, let $\varOmega$ be the set of points given by (\ref{eq9}) and let $|\varOmega|$ be its area.
		Then
		\begin{equation}
			|\varOmega| \leq (r/a)^{2/k}.
		\end{equation}
	\end{lemma}
	\noindent
	Now let $A$ is a $n\times n$ nonzero complex matrix having spectrum $\{\l_1, \dots, \l_n\}$, with $n\geq 3$ and suppose that we know a few distinct eigenvalues of $A$, namely $\l_1, \dots, \l_k$.
	Let $B$, $\rho(B)$, $r_j$, $r$, $D(\l_j,r_j)$ be exactly as in Theorem \ref{t1} and Theorem \ref{t2}.	
	According to theorem \ref{t1}, all the eigenvalues of $A$ are confined within the region of the complex plane given by
	\begin{equation}\label{eq11}
		\varOmega = \left\{z \in \mathbb{C}~~\big | ~~ \prod_{j=1}^{k} |z-\l_j| \leq \rho(B) \right\}.
	\end{equation}
	On the other hand, Lemma \ref{l3} implies that the size of this region satisfies the inequality
	\begin{equation}\label{eq8} 
		|\varOmega| \leq \pi (\rho(B))^{2/k}
	\end{equation}
	If $\rho(B) < r$, then it follows from (\ref{eq10}) and (\ref{eq8}) that
	$$ 
	|\varOmega| \leq \pi (\rho(B))^{2/k} < \pi (r)^{2/k} = \pi ~\sqrt[k]{\prod_{j=1}^{k} r_j^2}.
	$$
	By the AM-GM inequality (Lemma \ref{l2}), we have
	\begin{equation}
		\frac{r^2_1 + \dots + r^2_k}{k} ~~\geq~~ \sqrt[k]{\prod_{j=1}^{k} r_j^2}.
	\end{equation}
	Therefore,
	$$
	|\varOmega| ~~<~~ \frac{(\pi\, r_1^2) + \dots + (\pi\, r_k^2)}{k}.
	$$
	This complete the proof of the following theorem which is an extension of Theorem \ref{t2}.
	\begin{theorem}\label{t3}
		If $\rho(B) < r$, then
		$$
		|\varOmega| < \frac{(\pi\, r_1^2) + \dots + (\pi\, r_k^2)}{k},
		$$		
		that is, the spectrum of $A$ is confined within $k$ ovals, the sum of whose areas is $k$ times less than the  sum of the areas of the disks $D(\l_j,r_j), j=1, \dots, k$.
	\end{theorem}
	Theorem \ref{t3} has an important significance regarding the numerical computations of the eigenvalues $\l_1, \dots, \l_k$. 
	Suppose that $\s_1, \dots, \s_k$ are approximations of $\l_1, \dots, \l_k$ with absolute errors, respectively, $\e_1, \dots, e_k$.
	That is $|\l_i-\s_i|\leq \e_i$ for $i=1, \dots, k$.
	Consider the matrix
	\begin{equation}
		C = \prod_{j=1}^{k} (A - \s_j I).
	\end{equation}
	The eigenvalues of $C$ are given by $\prod_{j=1}^{k} (\l_i - \s_j)$ for $i=1, \dots, n$.
	For $\l_i \in \{\l_1, \dots, \l_k\}$, $\prod_{j=1}^{k} (\l_i - \s_j) \approx 0$ and the matrix $C$ reaches its spectral radius
	$\rho(C) = \prod_{j=1}^{k} |\l - \s_j|$ for some $\l$ in $\{\l_{k+1}, \dots, \l_n\}$.
	Then we have
	\begin{align}
		\rho(C)	& =  \prod_{j=1}^{k} |\l - \s_j|  \notag \\
				& \leq  \prod_{j=1}^{k} \left(|\l - \l_j| + \e_j \right) \notag\\
				& \leq  \prod_{j=1}^{k} \left(|\l - \l_j|\right) ~+~ \sum_{j=1}^{K} \left[\e_j \sum_{\substack{1\leq i \leq k\\i\neq j}} \l_i\right]
				~+~ \sum_{1\leq i < j \leq k} \mathcal{O} (\e_i \e_j). \label{17}\\ \notag
	\end{align}
	\noindent
	If  $\{\e_i\}$ are  small enough such that
	$$
	\sum_{1\leq i < j \leq k} \mathcal{O} (\e_i \e_j) ~~ \leq ~~ \sum_{j=1}^{K} \left[\e_j \sum_{\substack{1\leq i \leq k\\i\neq j}} \l_i\right]
	$$
	and
	$$
	\sum_{j=1}^{K} \left[\e_j \sum_{\substack{1\leq i \leq k\\i\neq j}} \l_i\right] ~~ << ~~ \prod_{j=1}^{k} \left(|\l - \l_j|\right),
	$$
	then it follows from (\ref{17}) that
	$$
	\rho(C) \approx  \prod_{j=1}^{k} |\l - \l_j|,
	$$
	that is,
	$$
	\rho(C) \approx \rho(B).
	$$
	Because of the numerical errors, some of the eigenvalues of $A$ might jump outside the region defined by
	\begin{equation}\label{eq12}
		\left\{ z \in \mathbb{C} ~~ | ~~ \prod_{j=1}^{k} |z-\sigma_j| \leq \rho(C) \right\},
	\end{equation}
	but in the case where $\rho(C) <r$, Theorem \ref{t3} gives confidence that  the spectrum of $A$ stays clustered within the disks obtained by using the  Euclidean minimum spanning tree of $\sigma_1, \dots, \sigma_k$ as the region given by (\ref{eq12}) is more than $k$ times smaller than the union of those disks and each oval is contained within one disk and centers around the same center.\\

\section{Larger matrices and the Lanczos algorithm}

We assume now that $A$ is hermitian or real symmetric.  Theorem~\ref{t1} can
provide interesting information about the distribution of eigenvalues, even
when a relatively small number of eigenvalues are known, as might be the case
when computing eigenvalues and eigenvectors using the Lanczos algorithm applied
to $A$.  If the number of eigenvalues used to form $B$ is relatively small, it is
also a simple matter to use the Lanczos algorithm to compute $\rho(B)$.  This
typically converges more quickly than power iteration and is a common way to
compute the spectral radius of a symmetric matrix.

To illustrate this point, we generated a $1000\times 1000$ diagonal matrix $D$ with
eigenvalues
\begin{equation*}
  \lambda_j = 0.5^{(j-500)/50}, \qquad \text{for $j=1, \ldots, 500$},
\end{equation*}
\begin{equation*}
  \lambda_j = 0.5^{(j-500)/50}+1000, \qquad \text{for $j=501, \ldots, 999$},
\end{equation*}
and $\lambda_{1000} = -1$.  This gives a relatively tight cluster of 500 eigenvalues in
$[0,1]$ and a more spread out cluster in $[1001, 2009.91]$ along with a more
isolated eigenvalue at $-1$.  We then formed a matrix $A=QDQ^T$, where $Q$ was
an orthogonal matrix randomly generated according to the Haar measure using the
approach of \cite{stew:80}.  Neither the eigenvalue regions nor the convergence
properties of the Lanczos algorithm depend on $Q$.  The eigenvalue $-1$ was
added to make sure the Lanczos algorithm could find an eigenvalue in or near
the cluster of eigenvalues in $[0,1]$ in a reasonable number of iterations.  The
tight clustering exhibited by the smaller eigenvalues makes them somewhat
difficult to find with Lanczos.  This would not be an issue if the smaller
eigenvalues were more spread out.

\begin{figure}
  \centering
  \vspace{1ex}
  \includegraphics[scale=0.30]{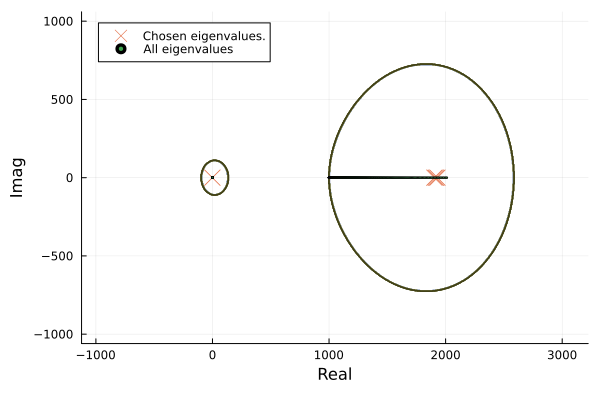}
  \caption{Eigenvalues and Theorem~\ref{t1} Regions.}
  \label{fig:lanczos_regions}
\end{figure}

\begin{figure}
  \centering
  \vspace{1ex}
  \includegraphics[scale=0.40]{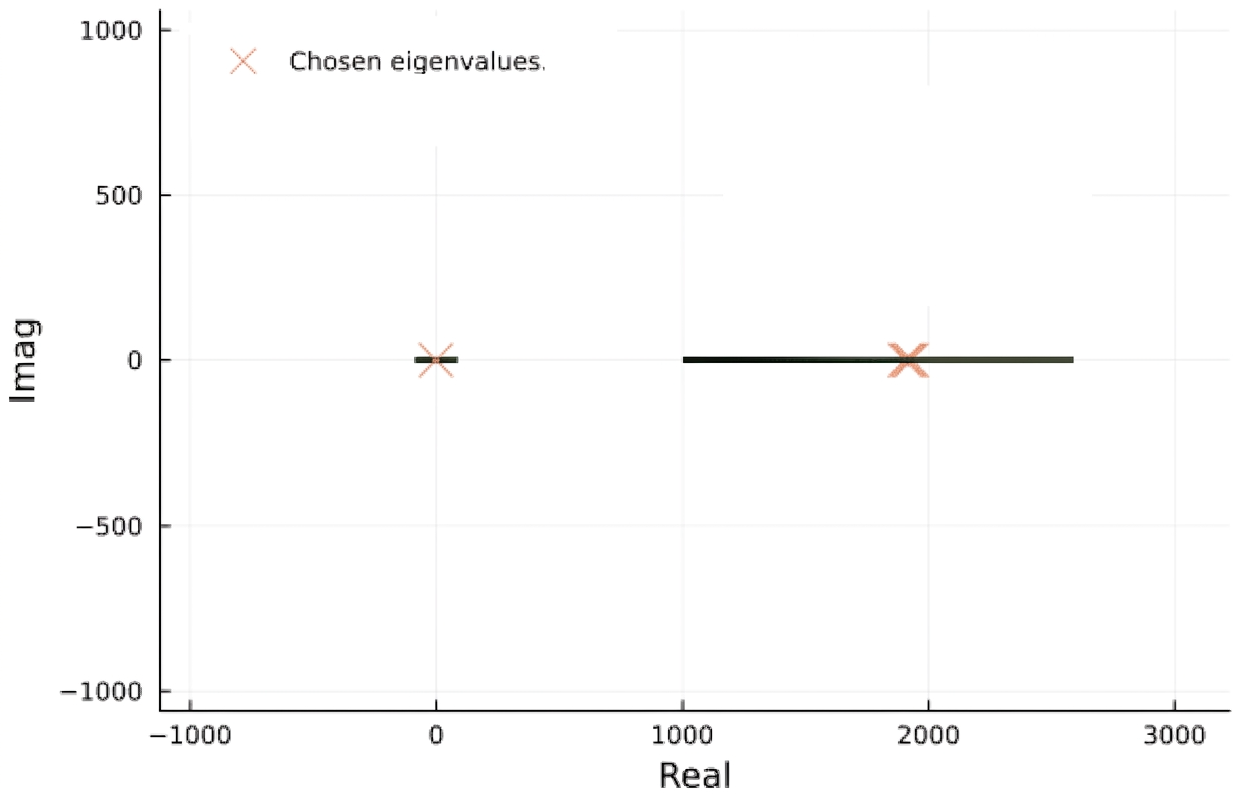}
  \caption{Selected Eigenvalues and Corollary~\ref{c2} Intervals}
  \label{fig:intervals}
\end{figure}

In 116 iterations the Lanczos algorithm obtained 11 eigenvalues with relative
residuals close to the unit round-off.  We selected $k=4$ in Theorem~\ref{t1}
and chose the four smallest identified eigenvalues
\begin{equation*}
  \lambda_1 = -0.9999999999988947, \quad
  \lambda_2 = 1891.4437768152313, \quad
\end{equation*}
\begin{equation*}
  \lambda_3 = 1903.8878682016305, \quad\text{and}\quad
  \lambda_4 = 1916.5056726302432.
\end{equation*}
The number of iterations was largely determined by the desire to find an
eigenvalue outside of the cluster of larger eigenvalues.  The other 7 of the 11
eigenvalues were the largest eigenvalues of $A$.  Note that these are all
relatively extreme eigenvalues.  We also used the Lanczos algorithm to compute
\begin{equation*}
  \rho(B) = 7.369855073399108\times 10^{11}.
\end{equation*}
The regions satisfying (\ref{eq5}) are shown in
Figure~\ref{fig:lanczos_regions}.  We have shown the regions in the complex
plane, but the eigenvalues contained in the regions are real.  Each eigenvalue
chosen to form $B$ is shown as a large ``$\times$.''  Note that, although the markers
overlap and are not entirely easy to see individually, there is one in the
smaller group of eigenvalues and 3 closely overlapping markers in the larger
group.  The other eigenvalues of $A$ are shown as small circles, which blend
together visually to form a thicker line segment in the large oval and a short
line segment in the smaller oval.  The intervals associated with Corollary~3.1,
are shown in Figure~\ref{fig:intervals}, with the darker lines marking the
intervals and the $\times$ markers showing the selected eigenvalues as before.
Given these selected eigenvalues, the intervals shown contain all of the other
eigenvalues.

The figures show that using four relatively extreme eigenvalues and the
spectral radius of $B$, the generalized Cassini ovals are able to resolve the
spectrum into two well separated clusters, revealing a large gap in the middle
of the spectrum.  While the results of this paper suggest that identifying
separated clusters of eigenvalues from a few computed eigenvalues should be
possible, one does not normally expect to get much information about
eigenvalues in the middle of the spectrum from the Lanczos algorithm.  From the
point of view of what the Lanczos algorithm usually provides, the fact that
Theorem~\ref{t1} and Corollary~\ref{c2} allow the identification of a gap in
the spectrum is a pleasant surprise.

In this section we have considered an example where the resulting region is
particularly nice and clearly illustrates clustering.  The results can vary
significantly with other matrices, but the approach has the advantages that the
regions are relatively easy to compute and are in many cases much smaller than
what could be found using Brauer sets or Gershgorin disks.  It is also possible
to do this with nonsymmetric problems and the Arnoldi iteration.  \vspace{12pt}

	\noindent
	\textbf{Acknowledgment.}
	The author Rachid Marsli would like to express sincere gratitude to Professor John D. King, Chair of the Department of Mathematics, Computer Science and Engineering at Georgia State University’s Perimeter College, for giving him the opportunity to return to academia and join the department as a faculty member in Fall 2026.


\begin{thebibliography}{30}
			\bibitem{A}					  Agarwal, P.K., Edelsbrunner, H., Schwarzkopf, O. et al. Euclidean minimum spanning trees and bichromatic closest pairs. Discrete Comput Geom
										  6 (1991), 407--422.
		
			\bibitem{Barany} 			 I. Barany and J. Solymosi, \textit{Gershgorin disks for multiple eigenvalues of non-negative matrices}, A Journey Through Discrete Mathematics,
										 Springer, Cham, 2017, 123--133.
			
			\bibitem{Brauer}		 	 A. Brauer, \textit{Limits for the characteristic roots of a matrix. II}, Duke Mathematical Journal, 14 (1947), 21--26.
			
								
			\bibitem{Dall}		 		 A. Dall'Acqua, D. Mugnolo, and M. Schelling, \textit{A New Gershgorin-Type Result for the Localisation of the Spectrum of Matrices},
										 Math. Nachr., 288 (2015) 1981--1994..
			
			\bibitem{De}		 		 G. V. Demidenko and Z. Wang, \textit{Localization of the Matrix Spectrum and Lyapunov Type Equations}, Lobachevskii J. Math., 45 (2024), 1808--1818.
			\bibitem{Chan}				 R. Chan, Circulant Preconditioners for Hermitian Toeplitz Systems, SIAM J. Matrix Anal. Appl. 10(1989), 542--550
			
			\bibitem{CS}			 	 R. Chan,  G. Strang, Toeplitz Equations by Conjugate Gradient with Circulant Preconditioner, SIAM J. Matrix Anal. Appl. 10(1989), 104--119
			
			\bibitem{Chu}				 E. K.-W. Chu and W.-W. Lin, Perturbation of eigenvalues for periodic matrix pairs via the Bauer–Fike theorems, Linear Algebra Appl. 
										 378 (2004), 183–202.
			
			\bibitem{G}	 		 		 S. A. Gerschgorin, \textit{\"{U}ber die Abgrenzung der Eigenwerte einer Matrix}, Izvestiya Rossiiskoi Akademii Nauk. Seriya
												 Matematicheskaya, (6) (1931), 749--754.
		
			\bibitem{Hall1}				 F. J. Hall and R. Marsli, \textit{On the location of eigenvalues of real constant row-sum matrices}, Bulletin of the Korean Mathematical
												 Society, 55 (2018), 1181--1187.
			
			\bibitem{Hall2}				 F. J. Hall, R. M. Marsli, and M. A. Stewart, \textit{An Application of Gelfand’s Formula in Approximating the Roots of Polynomials},
												 Experimental Mathematics, (2026), 1--13.		
												 
			\bibitem{Han}				 X. Han, X. Tong, Y. Fan, Eigen Selection in Spectral Clustering: A Theory-Guided Practice. Journal of the American Statistical
										 Association, 118(2023), 109–121. 
												 
			\bibitem{Horn}				 R. A. Horn and C. R. Johnson, Matrix Analysis, second edition, Cambridge University Press, Cambridge, 2013.
			
			\bibitem{K}	 				 L. Kolotilina, \textit{Generalizations of the Ostrowski–Brauer Theorem}, Linear Algebra Appl., 364 (2003), 65--80.
			
			\bibitem{L}					 U. von Luxburg, A tutorial on spectral clustering. Stat Comput 17, 395–416 (2007).
			
			\bibitem{M1}		 		 R. Marsli and F. J. Hall, \textit{Geometric Multiplicities and Ger\v{s}gorin Discs}, The American Mathematical Monthly, 120 (2013),
												 452--455.
			
			\bibitem{M2}				 R. Marsli and F. J. Hall, \textit{Inclusion regions and bounds for the eigenvalues of matrices with a known eigenpair}, Special Matrices,
												 8 (2020), 204--220.
			
			\bibitem{M3}				 R. Marsli and F. J. Hall, \textit{On the location of eigenvalues of real matrices}, Electronic Journal of Linear Algebra, 28 (2015),
												 206--213.
							
			\bibitem{M4} 				 R. Marsli and F. J. Hall, \textit{Some Refinements of Gersgorin Discs}, Int. J. Algebra, 7 (2013), 573--580.
		
			\bibitem{Me}				 A. Melman, \textit{Eigenvalue localization under partial spectral information}, Linear Algebra Appl., 573 (2019), 12--25.						
			
					
			\bibitem{Pena}				 J. M. Peña, \textit{Eigenvalue Localization for Symmetric Positive Toeplitz Matrices}, Axioms, 14 (2025), 232.	
			
			
			\bibitem{PR}				 F. P. Preparata and M. I. Shamos, Computational Geometry: An Introduction, Springer-Verlag, New York, 1985.
			
			\bibitem{Poyla}				 Pólya, G. (1928). Beitrag zur Verallgemeinerung des Verzerrungssatzes auf mehrfach zusammenhängende Gebiete. Sitzungsberichte der Preussischen
										 Akademie der Wissenschaften, Physikalisch-Mathematische Klasse, Berlin, 228–232 and 280–282.
			
			
			\bibitem{Rehman}			 M. Rehman, R. Amna, S. Akram, T. H. Rasulov, N. H. Mamatova, J. Sh. Ostonov, B. Z. Gafurov, and M. H. Boboqulova, 
												 \textit{Some New Localization Results for Eigenvalues of Complex Matrices}.
			
			\bibitem{S}				 	G. W. Stewart, Error and Perturbation Bounds for Subspaces Associated with Certain Eigenvalue Problems, SIAM Review, 15 (1973), 727–764.

                        \bibitem{stew:80} G. W. Stewart, The Efficient Generation of Random Orthogonal Matrices with an Application to Condition Estimators, SIAM J. Numer. Anal., 17(3):403--409, 1980.
			
			\bibitem{V}					 R. S. Varga, \textit{Gerschgorin and His Circles}, Springer Series in Computational Mathematics, Vol. 36, Springer-Verlag, Berlin, 2004.
			
			\bibitem{Tao}				 Tao Xiang and Shaogang Gong. Spectral clustering with eigenvector selection. Pattern Recognition, 41(3):1012–1029, 2008.
			
			\bibitem{W}					 J. Wu, D. Wang, and J. Zhao, \textit{Localization of Eigenvalues and Estimation of the Spread for Complex Matrices}, Izv. Math., 79 (2015),
												 208--215.
			
			\bibitem{Z}					 J. Zhao, \textit{A Survey of the Progress of Locating Methods of Complex Matrices' Eigenvalues and Some New Location Theorems and Their
												 Applications}, IMA J. Appl. Math., 80 (2014), 1273--1285.			
	\end{thebibliography}
\end{document}